\documentclass{amsart}
\usepackage{amssymb,amsmath,amsthm,graphics,amscd,amsfonts}

\theoremstyle{plain}
\newtheorem{theorem}{Theorem}[section]
\newtheorem{proposition}[theorem]{Proposition}
\newtheorem{corollary}[theorem]{Corollary}
\newtheorem{lemma}[theorem]{Lemma}

\newtheorem{definition}[theorem]{Definition}

\newtheorem{fact}[theorem]{Fact}
\newtheorem{remark}[theorem]{Remark}

\newcommand{\bfz}{{\bf z}}

\newcommand{\bfC}{{\mathbb C}}
\newcommand{\bfP}{{\mathbb P}}
\newcommand{\bfR}{{\mathbb R}}
\newcommand{\bfZ}{{\mathbb Z}}

\newcommand{\bfQ}{{\mathbb Q}}

\newcommand{\barg}{{\overline g}}

\newcommand{\barj}{{\overline j}}

\newcommand{\barz}{{\overline z}}

\newcommand{\barpartial}{{\overline \partial}}

\newcommand{\txi}{{\widetilde \xi}}

\newcommand{\mapright}[1]{\smash{\mathop{   \hbox to 0.7cm{\rightarrowfill}}
  \limits^{#1}}}

\def\bp{\overline{\partial}}

\def\CP{{\mathbb C \mathbb P}}

\def\p{\partial}
\def\bp{\overline{\partial}}

\def\vph{\varphi}

\def\la{\langle}
\def\ra{\rangle}

\begin{document}

\title{Momentum construction on Ricci-flat K\"ahler cones}
\author{Akito Futaki}
\address{Department of Mathematics, Tokyo Institute of Technology, 2-12-1,
O-okayama, Meguro, Tokyo 152-8551, Japan}
\email{futaki@math.titech.ac.jp}

\date{July 31, 2009 }

\begin{abstract} 
We extend Calabi ansatz over K\"ahler-Einstein manifolds  to Sasaki-Einstein manifolds.
As an application
we prove the existence of a complete scalar-flat K\"ahler metric on K\"ahler cone manifolds over
Sasaki-Einstein manifolds. In particular there exists a complete scalar-flat K\"ahler metric on
the toric K\"ahler cone manifold constructed from a toric diagram with 
a constant height. 
\end{abstract}



\maketitle

\section{Introduction}
A method to construct complete K\"ahler metrics with good curvature property is known as Calabi ansatz.
This method was employed first by Calabi \cite{calabi79} to construct complete Ricci-flat K\"ahler metrics on the canonical line
bundle over a K\"ahler-Einstein manifolds of positive scalar curvature. Calabi ansatz was extended by various
authors  (e.g. \cite{LeBrun88}, \cite{PedPoon91}, \cite{Simanca91}, \cite{Hwang-Singer02}). In this paper we extend Calabi ansatz to
the cone over Sasaki-Einstein manifolds. The method of this paper owes much to the work by Hwang and Singer (\cite{Hwang-Singer02})
who abstracted the essence of the earlier constructions on the total space of an 
Hermitian line bundle
$p : L \to M$ with ``$\sigma$-constant'' curvature over a K\"ahler manifold $(M, g_M)$. This method, also called
the {\it momentum construction}, searches for K\"ahler forms of the form
\begin{equation}\label{i1}
\omega = p^*\omega_M + dd^c f(t)
\end{equation}
where
$t$ is the logarithm of the norm function and $f$ is a function of one variable.
In this paper we extend
the momentum construction to Sasaki-Einstein manifolds or more generally 
Sasakian $\eta$-Einstein manifolds, i.e.  Sasaki
manifolds which are transversely K\"ahler-Einstein. 
Recall that a Sasaki manifold $(S,g)$ is an
odd dimensional  
Riemannian manifold such that its cone manifold 
$$(C(S), \barg) = (\bfR^+ \times S, dr^2 + r^2 g)$$
is a K\"ahler manifold where $r$ is the standard coordinate of $\bfR^+$. 
Then $S$, which is identified with the submanifold $\{r = 1\}$, becomes a contact manifold with the
contact form $\eta := (i(\barpartial - \p) \log r)|_{r=1} $. The vector field $\xi = J\frac \p{\p r}$ on $S$ is the
Reeb vector field of the contact form, that is
$$ i(\xi) \eta = 1 \quad \mathrm{and} \quad i(\xi) d\eta = 0.$$
We call the flow generated by $\xi$ the Reeb flow or the characteristic foliation. 
Since $\xi - iJ\xi$ is a holomorphic vector field and 
since $d\eta$ is non-degenerate on the kernel of $\eta$,
the local orbit spaces of the Reeb flow (or equivalently the local leaf spaces of the characteristic foliation) inherit K\"ahler structures.
It is a standard fact that $S$ has a Sasaki-Einstein metric 
if and only if the local orbit spaces of the Reeb flow have K\"ahler-Einstein metrics of
positive scalar curvature (positive K\"ahler-Einstein structure for short). If the K\"ahler form and
Ricci form on the local orbit spaces are denoted by $\omega^T$ and $\rho^T$ respectively
then the condition is expressed by
\begin{equation}\label{i2}
\rho^T = \kappa\, \omega^T
\end{equation}
for some positive constant $\kappa$. Allowing $\kappa$ to be any real constant a Sasaki manifold
with the condition
(\ref{i2}) is called an $\eta$-Einstein manifold, see section 2 for more detail.
Typical examples of Sasaki-Einstein manifolds are when $M$ is a positive
K\"ahler-Einstein manifold and the Sasaki manifold $S$ is the total space of 
the $U(1)$-bundle associated
with the canonical line bundle $K_M$. 
We perform the momentum construction on the cone $C(S)$ of a compact $\eta$-Einstein Sasaki manifold $S$
by replacing $\omega_M$ by the transverse K\"ahler form $\omega^T$ and 
putting $t = \log r$. Thus the momentum construction over an $\eta$-Einstein Sasaki manifold is of the form
\begin{equation}\label{i3}
\omega = p^*\omega^T + dd^c f(t)
\end{equation}
where $p : C(S) \to S$ is the projection along the radial flow generated by $r\frac \p{\p r}$.

Note that it is proved in \cite{FOW} (see also \cite{CFO} )
that given a compact toric Sasaki manifold $S$
corresponding to a toric diagram with a constant height, 
one can deform the Sasaki structure by varying the Reeb vector field 
and then deforming the transverse K\"ahler form to obtain a Sasaki-Einstein metric.
Given a toric Fano manifold $M$, 
we can apply this result to  the total space $S$ of the $U(1)$-bundle associated with the canonical line bundle $K_M$, and thus we obtain a Sasaki-Einstein metric on $S$ and a 
possibly non-standard vector Reeb field such that the local leaf spaces have positive K\"ahler-Einstein
metrics.

Typical results we obtain are as follows.

\begin{theorem}\label{main3}Let $S$ be a compact 
Sasaki-Einstein manifold and $C(S)$ its K\"ahler cone manifold with the cone metric $dr^2 + r^2g$.
Then we have the following.
\begin{enumerate}
\item[(a)]  There exists a complete scalar-flat K\"ahler metric on $C(S)$. 
\item[(b)] For any negative constant $c$ there exists a 
$\gamma > 0$ such that 
the submanifold $\{0 < r < \gamma \}$ in $C(S)$ 
admits a complete K\"ahler metric of negative constant scalar curvature $c$. 
\end{enumerate}
In particular if $C(S)$ is
a toric K\"ahler cone manifold corresponding to a toric diagram with a constant height then 
there exists a complete scalar-flat K\"ahler metric on $C(S)$.
\end{theorem}
\noindent
In this theorem no metric is Einstein.
As a special case of Theorem \ref{main3} we have the following corollary.

\begin{corollary}\label{main4}Let $M$ be a toric Fano manifold and $L$ be a holomorphic line bundle
such that $K_M = L^{\otimes p}$. Then for any positive integer $k$ the total space of $L^k$ minus the zero section 
admits a complete scalar-flat K\"ahler metric. There also exists a complete K\"ahler metric of negative constant scalar curvature on an 
open disk bundle of  $L^{\otimes k}$.
\end{corollary}

\begin{theorem}\label{main5}Let $S$ be compact $\eta$-Einstein Sasaki manifold 
with $\rho^T = \kappa \omega^T$ for some non-positive constant $\kappa$
and $C(S)$ its cone 
manifold with cone metric $dr^2 + r^2g$. 
\begin{enumerate}
\item[(a)]  
If $\kappa = 0$, then for any negative constant $c$, there exists a 
$\gamma > 0$ such that the submanifold $\{0 < r < \gamma \}$ in $C(S)$ 
admits a complete K\"ahler metric of negative constant scalar curvature $c$. 
\item[(b)]  
If $\kappa < 0$ we have a negative constant $c_0$ such that
there exists a complete K\"ahler metric on $C(S)$ with scalar curvature $c_0$ and that
for any negative constant $c < c_0$ there exists a 
$\gamma > 0$ such that the submanifold $\{0 < r < \gamma \}$ in $C(S)$ 
admits a complete K\"ahler metric of negative constant scalar curvature $c$. 
This metric is Einstein if $c = (m+1)\kappa$.
\end{enumerate}
\end{theorem}
\noindent
Note that in the case of (a), there is no borderline case and no metric is Einstein, 
as the proof given in section 4 shows.
There are many examples of compact Sasaki manifolds with
$\rho^T = \kappa \omega^T$ for some non-positive constant $\kappa$, 
see Remark \ref{main7} in section 4 or Boyer, Galicki and Matzeu \cite{BGM} for more detail. 
Theorem \ref{main5} therefore produces many new complete scalar-flat K\"ahler
manifolds and negative K\"ahler-Einstein manifolds.

As in the case of Calabi \cite{calabi79}, one may try to 
prove the existence of a complete Ricci-flat 
K\"ahler metric on the total space of the canonical line bundle $K_M$
of a toric Fano manifold $M$.
In fact one can prove the following result.

\begin{proposition}\label{main2}Let $M$ be a toric Fano manifold and $L$ be
a holomorphic line bundle over $M$ such that $K_M = L^{\otimes p}$ for some positive integer $p$.
Then for each positive integer $k$, there exists a scalar-flat K\"ahler metric
on the total space of $L^{\otimes k} - \mathrm{zero\ section}$ . This metric is Ricci-flat when $k = p$, that is
when $L^{\otimes k} = K_M$.
\end{proposition}
It is not clear if the metric extends smoothly to the zero section if the Sasaki-Einstein structure on the $U(1)$-bundle of $K_M$ is
irregular. We will give some consideration on this point in section 4.
 In physics literature Oota and Yasui \cite{OY06} have obtained a complete Ricci-flat K\"ahler
 metric on the canonical bundle of one-point-blow-up of $\bfC\bfP^2$ by a different
 derivation.

The organization of this paper is as follows. In section 2 we review Sasakian geometry and
give precise statements of the results obtained in \cite{FOW}. In section 3 we apply the Calabi
ansatz to Sasakian $\eta$-Einstein manifolds and derive basic formulae. In section 4 we give proofs of
the theorems stated in this introduction. In section 5 we give a proof of Proposition \ref{main2} and other related results, and also discuss about
the possibility of extending the metric to the zero section.

We refer the reader to the review of Boyer and Galicki \cite{BG07} in which the problem
of resolving the cone metrics are taken up from wider view points.

This work was inspired by conversations with Shi-Shyr Roan and Naichung Conan Leung at
Sugadaira conference in 2006. The author is especially grateful to Prof. Leung 
for suggesting to consider the canonical bundles of toric Fano manifolds.

\section{Sasaki manifolds} 

A Sasaki manifold is a Riemannian manifold $(S, g)$ whose cone manifold 
$(C(S), \overline g)$ with $C(S) \cong  S\times \bfR^+$ and $\overline g = dr^2  + r^2g$ 
is a K\"ahler manifold where $r$ is the standard coordinate on $\bfR^+$.
From this definition $S$ is odd-dimensional and we put $\dim _{\bfR}S = 2m + 1$, and thus
$\dim_{\bfC} C(S) = m+1$.
Note that $C(S)$ does not contain the apex. $S$ is always identified with the
real hypersurface $\{r=1\}$ in $C(S)$.

Putting $\txi = J(r\frac{\p}{\partial r})$, 
$\txi - iJ\txi$ defines a holomorphic vector field on $C(S)$. 
The restriction of $\txi$ to 
$S \cong \{r=1\}$, which is tangent to $S$, is called the {\it Reeb vector field} of $S$ and
denoted by $\xi$. 
The Reeb vector field $\xi$ is a Killing vector field on $S$, and thus generates a toral subgroup 
in the isometry group of $S$.

Let $\eta$ be the dual $1$-form
to $\xi$ using the Riemannian metric $g$. Then $\eta$ can be expressed as
$$\eta = (i(\barpartial - \p) \log r)|_{r=1} =  (2d^c \log r) |_{r=1}.$$ 
Put $D = \mathrm{Ker}\, \eta$ and call it the {\it contact bundle}.
Then $d\eta$ is non-degenerate on $D$ and thus $S$ becomes a contact manifold
with the contact form $\eta$.
The Reeb vector field $\xi$ satisfies
$$i(\xi)\eta= 1\quad \mathrm{and}\quad i(\xi) d\eta = 0$$
where $i(\xi)$ denotes the inner product, which are often used as the defining
properties of the Reeb vector field for contact manifolds.

The Reeb vector field $\xi$ generates a $1$-dimensional foliation $\mathcal F_{\xi}$, called the characteristic 
foliation, on $S$. We also regard $\mathcal F_{\xi}$ as the flow generated by $\xi$ and call it the Reeb flow. 
Since 
the Reeb flow shares common local orbit spaces with
the holomorphic flow generated by $\txi - iJ\txi$ on $C(S)$
and the latter has a natural transverse holomorphic structure
the characteritic foliation on $S$ admits  
a transverse holomorphic structure. The tangent spaces to the local leaf spaces are naturally
isomorphic to a fiber of $D$, and using this isomorphism and considering the symplectic form
$\frac 12 d\eta$ on $D$ we obtain a well-defined K\"ahler form on the local leaf spaces of 
$\mathcal F_{\xi}$.
Though the K\"ahler forms are defined on the local leaf spaces they are pulled back to $S$ and
glued together to give a global $2$-form 
\begin{equation}\label{s1}
\omega^T = \frac12 d\eta = d(d^c \log r\,|_{r=1}) = (dd^c \log r)|_{r=1}
\end{equation}
on $S$, which we call the {\it transverse K\"ahler form}. Note that the clumsy constant
$1/2$ is necessary since
\begin{equation}\label{s2}
d\eta (X,Y) = 2g(\Phi X,Y)
\end{equation}
for $X,\ Y \in D_x,\ x\in S$, where $\Phi$ is the natural complex structure on $D$.
We call the collection of K\"ahler structures on local leaf spaces of $\mathcal F_{\xi}$
the {\it transverse K\"ahler structure}. 

Recall that a smooth differential form $\alpha$ on $S$ is basic if
$$ i(\xi)\alpha = 0\quad \mathrm{and}\quad \mathcal L_{\xi}\alpha = 0$$
where $\mathcal L_{\xi}$ denotes the Lie derivative by $\xi$.
The basic forms are preserved by the exterior derivative $d$ which decomposes into
$d = \p_B + \barpartial_B$, and we can define basic cohomology groups and
basic Dolbeault cohomology groups. We also have the transverse Chern-Weil theory and
can define basic Chern classes for complex vector bundles with basic transition functions.
As in the K\"ahler case the basic first Chern class
of the Reeb foliation $c_1^B$ is represented by the $1/2\pi$ times the transverse
Ricci form $\rho^T$:
\begin{equation}
\rho^T = -i \p_B\bp_B \log \det(g^T_{i\barj})
\end{equation} 
where 
$$\omega^T = i \ g_{i\barj}^T\ dz^i \wedge dz^{\barj}$$
and $z^1,\ \cdots,\ z^m$ are local holomorphic coordinates on the
local leaf space.

Now we turn to the study of Sasaki-Einstein manifolds. We start with
the following fact.
\begin{fact}\label{s3}Let $(S,g)$ be a $(2m+1)$-dimensional Sasaki manifold. 
The following three conditions are equivalent.
\begin{enumerate}
\item[(a)] $(S,g)$ is a Sasaki-Einstein manifold. The Einstein constant is necessarily $2m$.
\item[(b)] $(C(S), \overline g)$ is a Ricci-flat K\"ahler manifold.
\item[(c)] The local leaf spaces of the Reeb foliation have transverse K\"ahler-Einstein metrics
with Einstein constant $2m+2$.
\end{enumerate}
\end{fact}
\noindent
See for proofs \cite{boga99} or \cite{BGbook}. A typical example of Fact \ref{s3} is when 
$$(C(S), S, \mathrm{local\ leaf\ spaces}) = (\bfC^{m+1}-\{0\}, S^{2m+1}, \bfC\bfP^m)$$ 
where $S^{2m+1}$ is the standard sphere of dimension $2m+1$. Of course the Reeb flow is induced by
the standard $S^1$-action.

Suppose that $S$ has an Einstein metric. Then by (c) of Fact \ref{s3} we have
\begin{equation}\label{s4}
\rho^T = (2m+2) \omega^T = (m+1)d\eta,
\end{equation}
hence $c_1^B > 0$, i.e. $c_1^B$ is represented by a positive basic $(1,1)$-form. Moreover
under the natural homomorphism $H^2_B(\mathcal F_{\xi}) \to H^2(S)$ of the basic cohomology
group to the ordinary de Rham cohomology group the basic first Chern class $c_1^B$ is sent
to the ordinary first Chern class $c_1(D)$, but by (\ref{s4}) 
\begin{equation}\label{s5}
c_1(D) = (2m+2)\omega^T = (m+1)[d\eta] = 0.
\end{equation}
Conversely if $c_1^B > 0$ and $c_1(D) = 0$ then $c_1^B = \tau [d\eta]$ for some positive
constant $\tau$. See Proposition 4.3 in \cite{FOW} for the detail.

A Sasaki manifold 
$(S, g)$ is said to be toric if the K\"ahler cone manifold $(C(S), \barg)$ is toric, namely if
$(m+1)$-dimensional torus $G$ acts on $(C(S), \overline g)$ effectively as
holomorphic isometries. Note that then $G$ preserves $\txi$ because $G$ preserves $r$ and the
complex structure $J$. Since the one-parameter group of transformations generated by 
$\txi$ acts on $C(S)$ as holomorphic isometries
and since $G$ already has the maximal dimension of possible torus actions on $C(S)$
then $\txi$ is contained in the Lie algebra $\mathfrak g$ of $G$. The action of $G$ on $C(S)$
naturally descends to an action on $S$.

\begin{definition}\label{good} Let $\mathfrak g^{\ast}$ be the dual of the Lie algebra $\mathfrak g$
of the $(m+1)$ dimensional torus $G$. Let $\bfZ_{\mathfrak g}$ be the integral lattice of $\mathfrak g$, that is
the kernel of the exponential map $\exp : \mathfrak g \to G$. 
A subset $C \subset \mathfrak g^{\ast}$ is a rational polyhedral cone if there exists a
finite set of vectors $\lambda_i \in \bfZ_{\mathfrak g}$, $1 \le i \le d$, such that
$$ C = \{ y \in \mathfrak g^{\ast}\ |\ \la y, \lambda_i \ra \ge 0\ \mathrm{for\ }\ i = 1, \cdots, d\}.$$
We assume that the set $\lambda_i$ is minimal in that for any $j$
$$ C \ne \{ y \in \mathfrak g^{\ast}\ |\ \la y, \lambda_i \ra \ge 0\ \mathrm{for\ all}\ i\ne j\}$$
and that each $\lambda_i$ is primitive, i.e. $\lambda_i$ is not of the form $\lambda_i = a\mu$
for an integer $a \ge 2$ and $\mu \in \bfZ_{\mathfrak g}$. 
(Thus $d$ is
the number of codimension $1$ faces if $C$ has non-empty interior.)
Under these two assumptions
a rational polyhedral cone $C$ with nonempty interior is good if the following condition holds.
If
$$ \{ y \in C\ |\ \la y, \lambda_{i_j}\ra = 0\ \mathrm{for\ all}\ j = 1, \cdots, k \}$$
is a non-empty face of $C$ for some $\{i_1, \cdots, i_k\} \subset \{1, \cdots, d\}$, then 
$\lambda_{i_1}, \cdots, \lambda_{i_k}$ are linearly independent
over $\bfZ$ and 
\begin{equation}\label{goodcondition}
\{ \sum_{j=1}^k a_j \lambda_{i_j}\ |\ a_j \in \bfR\} \cap \bfZ_{\mathfrak g} = 
\{ \sum_{j=1}^k m_j \lambda_{i_j}\ |\ m_j \in \bfZ\}.
\end{equation}
\end{definition}

Given a good rational polyhedral cone $C$ we can construct a {\it smooth} toric contact manifold
whose moment map image is $C$.
\begin{definition}\label{TD1} An $(m+1)$-dimensional toric diagram with height $\ell$ is a
collection of $\lambda_i \in \bfZ^{m+1} \cong 
\bfZ_{\mathfrak g}$ satisfying $(\ref{goodcondition})$ and 
$\gamma  \in \bfQ^{m+1} \cong
(\bfQ_{\mathfrak g})^{\ast}$
 such that 
\begin{enumerate}
\item[(1)] $\ell$ is a positive integer such that $\ell\gamma$ is a primitive element of
the integer lattice $\bfZ^{m+1} \cong \bfZ^{\ast}_{\mathfrak g}$.
\item[(2)] $\la \gamma, \lambda_i\ra = -1$. 
\end{enumerate}
We say that a good rational polyhedral
cone $C$ is associated with a toric diagram of height $\ell$ if there exists a rational vector $\gamma$ 
satisfying $(1)$ and $(2)$ above.
\end{definition}

The reason why we use the terminology ``height $\ell$'' is because 
using a transformation by an element of $SL(m+1,\bfZ)$ we may assume that
$$ \gamma = \left(\begin{array}{r} -\frac1\ell \\ 0 \\ \vdots \\ 0 \end{array}\right)$$
and the first component of $\lambda_i$ is equal to $\ell$ for each $i$.
The following theorem asserts that the condition of constant height in Definition \ref{TD1}
is a Calabi-Yau condition for $C(S)$.

\begin{theorem}[\cite{CFO}]\label{s6} Let $S$ be a compact toric Sasaki manifold with 
$\dim S \ge 5$. 
Then the following three conditions are equivalent.
\begin{enumerate}
\item[(a)] $c_1^B > 0$ and $c_1(D) = 0$.
\item[(b)] The Sasaki manifold $S$ is obtained from a toric diagram with height $\ell$ for some
positive integer $\ell$ defined by $\lambda_1, \cdots, \lambda_d \in \mathfrak g$ and
$\gamma \in \mathfrak g^{\ast}$ and the Reeb field
$\xi \in \mathfrak g$ satisfies
$$ \la \gamma, \xi \ra = -m-1\ \ \mathrm{and}\ \ \la y, \xi \ra > 0\ \ \mathrm{for\ all}\ 
y \in C$$
where $C = \{y \in \mathfrak g^{\ast}| \la y, \lambda_j \ra \ge 0,\ j=1,\cdots, d\}$.
\item[(c)]  For some positive integer $\ell$, the $\ell$-th power 
$K^{\otimes\ell}_{C(S)}$ of the canonical line bundle $K_{C(S)}$
 is trivial.
\end{enumerate}
\end{theorem}
The main theorem of \cite{FOW} is the following.

\begin{theorem}[\cite{FOW}]\label{s7} Suppose that we are
given a compact toric Sasaki manifold satisfying one of the equivalent conditions in 
Theorem \ref{s6} so that we may assume that $c_1^B = (2m+2)[\omega^T]$ as basic 
cohomology classes.
Then we can deform the 
Sasaki structure by varying the Reeb vector field and then performing a transverse
K\"ahler deformation to obtain a Sasaki-Einstein metric.
\end{theorem}

The proof of Theorem \ref{s7} is outlined as follows.  
Fixing a Reeb vector field we have a fixed transverse holomorphic structure. 
The first step is to prove that for the fixed transverse holomorphic structure we
can deform the transverse K\"ahler structure in the form
\begin{equation}\label{s10'}
\omega^T + dd^c \psi_1 = d(d^c (\log r + \psi_1)|_{r=1}) = \frac 12 d((d^c \log (r^2 \exp(2\psi_1))|_{r=1})
\end{equation}
where $\psi_1$ is a smooth basic function $S$ such that the deformed transverse K\"ahler metric satisfies the K\"ahler Ricci soliton equation:
$$ \rho^T  - (2m+2)\omega^T = \mathcal L_X \omega $$
for some ``Hamiltonian-holomorphic vector field'' $X$ in the sense of \cite{FOW}. 
Note that this deformation deforms $r^2$ into $r^2 \exp(2\psi_1)$.
A transverse K\"ahler-Ricci
soliton becomes a positive transverse K\"ahler-Einstein metric, i.e. $X = 0$, if the invariant $f_1$ which obstructs
the existence of positive transverse K\"ahler-Einstein metric vanishes. 
Note that $f_1$ depends only on the Reeb vector field $\xi$.
The second step is then to show that there exists a unique Reeb vector field $\xi'$
such that $f_1$ vanishes; this idea is due to Martelli-Sparks-Yau \cite{MSY2}. For this
$\xi'$ we have a positive transverse K\"ahler-Einstein metric, that is a Sasaki-Einstein metric.

\begin{definition}\label{eta}
A Sasaki metric $g$ is said to be $\eta$-Einstein if there exist constants $\lambda$ and
$\nu$ such that 
$$ \mathrm{Ric}_g = \lambda\, g + \nu\, \eta\otimes \eta$$
where $ \mathrm{Ric}_g$ denotes the Ricci curvature of $g$. 
\end{definition}
By elementary computations in Sasakian geometry 
we always have $\mathrm{Ric}_g(\xi,\xi) = 2m$ on any Sasaki manifolds.
This implies that
$\lambda + \nu = 2m$.
Let $\mathrm{Ric^T}$ denote the Ricci curvature of the local leaf space. Then again elementary
computations show
\begin{equation}\label{s10}
\mathrm{Ric}_g(\widetilde X,\widetilde Y) = (\mathrm{Ric^T} - 2g^T)(X,Y)
\end{equation}
where $\widetilde X,\ \widetilde Y \in D$ are the lifts of tangent vectors  $X,\ Y $ of local leaf space. 
From this we see that the condition of $\eta$-Einstein metric is equivalent to
\begin{equation}\label{s11}
\mathrm{Ric^T} = (\lambda + 2)g^T.
\end{equation}

Given a Sasaki manifold with the K\"ahler cone metric $\barg = dr^2 + r^2 g$, we transform
the Sasakian structure by deforming $r$ into $r' = r^a$ for positive constant $a$. 
This transformation is called the $D$-homothetic transformation. Then the new
Sasaki structure has 
\begin{equation}\label{s12}
\eta' = d \log r^a = a\eta, \quad \xi' = \frac 1a \xi,
\end{equation}
\begin{equation}\label{s13}
g' = a g^T + a\eta \otimes a\eta = ag + (a^2 -a)\eta\otimes\eta.
\end{equation}
Suppose that $g$ is $\eta$-Einstein with $\mathrm{Ric}_g = \lambda g + \nu \eta \otimes \eta$. 
Since the Ricci curvature of a K\"ahler manifold is invariant under homotheties
we have $\mathrm{Ric}^{\prime T} = \mathrm{Ric}^T$. From this and $\mathrm{Ric}_{g'}(\xi',\xi') = 2m$ we have
\begin{eqnarray}\label{s14}
\mathrm{Ric}_{g'} &=& \mathrm{Ric}^{\prime T} - 2g^{\prime T} + 2m \eta' \otimes \eta' \\
&=& \mathrm{Ric}^T - 2a g^T + 2m \eta' \otimes \eta' \nonumber\\
&=& \mathrm{Ric}|_{D \times D} + 2g^T - 2a g^T + 2m \eta' \otimes \eta' \nonumber\\
&=& \lambda g^T  + 2g^T - 2a g^T + 2m \eta' \otimes \eta'. \nonumber
\end{eqnarray}
This shows that $g'$ is $\eta$-Einstein with
\begin{equation}\label{s15}
\lambda' =  \frac{\lambda + 2 - 2a}a.
\end{equation}
In summary, we have
\begin{itemize}
\item Sasaki-Einstein metric is a special case of an $\eta$-Einstein metric with $\lambda = 2m$ and
$\nu = 0$;
\item a Sasaki metric is an $\eta$-Einstein metric with $\lambda + 2 > 0$ if and only if its transverse
K\"ahler metric is positive K\"ahler-Einstein with Einstein constant $\lambda + 2$;
\item
under the $D$-homothetic transformation of an $\eta$-Einstein metric we have
a new $\eta$-Einstein metric with
\begin{equation}\label{s16}
\rho^{\prime T} = \rho^T, \quad \omega^{\prime T} = a\omega^T, \quad
\rho^{\prime T} = (\lambda' + 2) \omega^{\prime T} = \frac{\lambda + 2}a \omega^{\prime T},
\end{equation}
and thus, for any positive constants $\kappa$ and $\kappa'$, a transverse K\"ahler-Einstein metric
with Einstein constant $\kappa$ can be transformed by a $D$-homothetic transformation to a
transverse K\"ahler-Einstein metric with Einstein constant $\kappa'$.
\end{itemize}

\section{Calabi ansatz for Sasakian $\eta$-Einstein manifolds}

Let $(S, g)$ be a Sasakian $\eta$-Einstein manifold with 
$\mathrm{Ric}_g = \lambda\,g + \nu\,\eta \otimes \eta$ 
and with K\"ahler cone metric on $C(S)$
$$ \barg = dr^2 + r^2g.$$
Let $\omega^T = \frac12 d\eta$ be the transverse
K\"ahler form which gives positive K\"ahler-Einstein metrics on local leaf spaces with
\begin{equation}\label{C0}
 \rho^T = \kappa \omega^T
 \end{equation}
where we have set 
$$\kappa := \lambda + 2.$$
As we work on $C(S)$ it is convenient to lift $\eta$ on $S$ to $C(S)$ by
\begin{equation}\label{C1}
\eta = 2 d^c \log r.
\end{equation}
We use the same notation $\eta$ for this lifted one to $C(S)$. Then $\omega^T$ is also lifted to
$C(S)$ by
\begin{equation}\label{C2}
\omega^T = dd^c \log r.
\end{equation}
Again we use the same notation $\omega^T$ for the lifted one to $C(S)$. The {\it Calabi ansatz}
searches for a K\"ahler form on $C(S)$ of the form
\begin{equation}\label{C3}
\omega = \omega^T + i \p\bp\, F(t)
\end{equation}
where $t = \log r$ and $F$ is a smooth function of one variable on $(t_1, t_2) \subset (-\infty, \infty)$.

We set
\begin{equation}\label{C4}
\tau = F'(t),
\end{equation}
\begin{equation}\label{C5}
\varphi(\tau) = F''(t).
\end{equation}
Since we require $\omega$ to be a positive form and
\begin{eqnarray}\label{C6}
i\p\bp\, F(t) &=& i\,F''(t)\, \p t \wedge \bp t \,+ i\, F'(t)\, \p\bp t \\
&=& i\,\varphi(\tau)\, \p t\wedge \bp t \,+ \tau\, \omega^T. \nonumber
\end{eqnarray}
then we must have $\varphi(\tau) > 0$.
We also require that the image of $F'$ is an open interval $(0, b)$ with $b \le \infty$, i.e.
\begin{equation}\label{C9}
\lim_{t \to t_1} F'(t) = 0, \qquad \lim_{t\to t_2} F'(t) = b.
\end{equation}
It follows from $\varphi(\tau) > 0$ that $F'$ is a diffeomorphism from $(t_1,t_2)$ to $(0,b)$.

\begin{definition}\label{C7}
We call $\varphi(\tau)$ the profile of the Calabi ansatz (\ref{C3}).
\end{definition}

Conversely, given a positive function $\varphi > 0$ on $(0,b)$ 
such that 
\begin{equation}\label{C11}
\lim_{\tau \to 0^+} \int_{\tau_0}^{\tau} \frac{dx}{\varphi(x)} = t_1, \qquad
\lim_{\tau \to b^-} \int_{\tau_0}^{b} \frac{dx}{\varphi(x)} = t_2
\end{equation}
we can recover the Calabi ansatz
as follows. 
Fix $\tau_0$ and introduce a function $\tau(t)$ by
\begin{equation}\label{C8}
t = \int_{\tau_0}^{\tau(t)} \frac{dx}{\varphi(x)},
\end{equation}
and then $F(t)$ by
\begin{equation}\label{C10}
F(t) = \int_{\tau_0}^{\tau(t)} \frac{xdx}{\varphi(x)}.
\end{equation}
Putting $t = \log r$ we may regard $\tau$ and $F$ as functions on
\begin{equation}\label{C12}
C(S)_{(t_1,t_2)} := \{ e^{t_1} < r < e^{t_2} \} \subset C(S),
\end{equation}
Put 
\begin{eqnarray}\label{C13}
\omega_{\varphi} &:=& \omega^T + dd^c\,F(t) \\
&=& (1+ \tau)\, \omega^T + \varphi(\tau)\,i \p t\wedge \bp t \nonumber \\
&=& (1+ \tau)\, \omega^T + \varphi(\tau)^{-1}\,i \p \tau\wedge \bp \tau\nonumber
\end{eqnarray}
As we assume $\varphi > 0$ on $(0,b)$, $\omega_{\varphi}$ defines a K\"ahler form
and have recovered Calabi ansatz.

\begin{remark}
When $S$ is the total space of the $U(1)$-bundle associated with an Hermitian line bundle $(L,h) \to M$
such that $i\p\bp \log h$ is a K\"ahler form on $M$, then we may regard the total space of $L$ 
as a Hamiltonian $U(1)$-space with $\tau = F'(t)$ the moment map. The K\"ahler potential $F$ along
the fiber is transformed under the Legendre transform into the symplectic potential $G$ given by the
symmetrical relation
\begin{equation}\label{C14}
G(\tau) + F(t) = \tau t.
\end{equation}
It is easy to see that $G''(\tau) = 1/\varphi(\tau)$. It is well-known that,
especially in the case of toric K\"ahler manifolds (c.f. \cite{Abreu}), the scalar curvature has
a less complicated expression if one uses the symplectic potential. 
This principle works in our situation as is  seen in the computations below.
\end{remark}

Next we compute the Ricci form $\rho_{\varphi}$ and the scalar curvature $\sigma_{\varphi}$ of 
$\omega_{\varphi}$. First we need to choose local holomorphic coordinates $(z^0, z^1, \cdots, z^m)$
on $C(S)$.
We choose $z^0$ to be the coordinate along the holomorphic Reeb flow, more precisely
\begin{equation}\label{C15}
\frac \p{\p z^0} = \frac 12 (r\frac \p{\p r} - iJ(r\frac \p{\p r})) = \frac12  (r\frac \p{\p r} - i\txi),
\end{equation}
and $z^1, \cdots, z^m$ to be the pull-back of local holomorphic coordinates on the local leaf space.
Then it is easy to check that
\begin{equation}\label{C16}
dz^0 = \frac{dr}r + i\eta,
\end{equation}
and that
\begin{equation}\label{C17}
idz^0 \wedge d\barz^0 = 2 \frac{dr}r \wedge \eta.
\end{equation}

Using these coordinates one can compute the volume form as
\begin{eqnarray}\label{C18}
\omega_{\varphi}^{m+1} &=& (1+\tau)^m(m+1) \varphi(\tau)\,i\p t \wedge\bp t \wedge(\omega^T)^m\\
&=&  (1+\tau)^m (m+1) \varphi(\tau)\,dt\wedge d^ct\wedge(\omega^T)^m\nonumber\\
&=&  (1+\tau)^m (m+1) \varphi(\tau)\,i\frac{dr}r\wedge \eta\wedge(\omega^T)^m\nonumber\\
&=&  (1+\tau)^m (m+1) \varphi(\tau)\,\frac i2 dz^0\wedge d\barz^0\wedge(\omega^T)^m\nonumber.
\end{eqnarray}
The Ricci form and the scalar curvature can be computed as
\begin{eqnarray}\label{C18'}
\rho_{\varphi} &=& \rho^T - i\p\bp \log ((1+\tau)^m \varphi(\tau)) \\
&=& \kappa \omega^T - i\p\bp \log ((1+\tau)^m \varphi(\tau)),\nonumber
\end{eqnarray}
\begin{eqnarray}\label{C19}
\sigma_{\varphi} &=& \frac{\sigma^T}{1+\tau} - \Delta_{\varphi} \log ((1+\tau)^m \varphi(\tau))\\
&=& \frac{m\kappa}{1+\tau} - \Delta_{\varphi} \log ((1+\tau)^m \varphi(\tau)),\nonumber
\end{eqnarray}
where $\Delta_{\varphi}$ denotes the complex Laplacian with respect to $\omega_{\varphi}$.

Let $u(\tau)$ be a smooth function of $\tau$. Then
\begin{eqnarray}\label{C20}
dd^c\, u(\tau) &=& d\,(u'(\tau)\, \frac{d\tau}{dt}\, d^c t)\\
&=& u'(\tau) \varphi(\tau)dd^c t + (u'\varphi)' \varphi dt \wedge d^ct \nonumber\\
&=& u'(\tau) \varphi(\tau)dd^c t + \frac 1{\varphi}(u'\varphi)'  d\tau \wedge d^c\tau. \nonumber
\end{eqnarray}
Taking wedge product of this with 
\begin{equation}\label{C21}
\omega_{\varphi}^m = (1+\tau)^m(\omega^T)^m + m(1+\tau)^{m-1}\wedge \varphi^{-1}d\tau
\wedge d^c\tau
\end{equation}
and comparing it with 
\begin{equation}\label{C22}
\omega_{\varphi}^{m+1} = (1+\tau)^m(m+1) \varphi^{-1}\,d \tau 
\wedge d^c \tau \wedge(\omega^T)^m.
\end{equation}
we obtain
\begin{equation}\label{C23}
\Delta_{\varphi} u = \frac m{1+\tau} u' \varphi + (u'\varphi )'. 
\end{equation}

Apply (\ref{C23}) with $u = \log ((1+\tau)^m \varphi)$ and insert it to (\ref{C19}). Then we obtain
\begin{eqnarray}\label{C24}
\sigma_{\varphi} &=& \frac{m\kappa}{1+\tau} - \Delta_{\varphi} \log ((1+\tau)^m \varphi)\\
&=& \frac{m\kappa}{1+\tau} - \frac m{1+\tau}\varphi\frac d{d\tau} \log ((1+\tau)^m \varphi)\\
&& \qquad \qquad\qquad - \frac d{d\tau}(\varphi \frac d{d\tau} \log (1+\tau)^m\varphi(\tau)) \nonumber\\
&=& \frac{m\kappa}{1+\tau} - \frac 1{(1+\tau)^m}\frac d{d\tau} ((1+\tau)^m\varphi \frac d{d\tau}
\log ((1+\tau)^m \varphi) \\
&=& \frac{m\kappa}{1+\tau} - \frac 1{(1+\tau)^m}\frac {d^2}{d\tau^2} ((1+\tau)^m\varphi).\nonumber
\end{eqnarray}
Setting $\sigma_{\varphi} = c$ we get an ordinary differential equation
\begin{equation}\label{C25}
(\varphi (1+\tau)^m)'' = \left(\frac{m\kappa}{(1+\tau)} - c\right)(1+\tau)^m.
\end{equation}
We can easily solve this equation as
\begin{equation}\label{C26}
\varphi(\tau) = \frac{\kappa}{m+1}(1+\tau) - \frac c{(m+1)(m+2)}(1+\tau)^2 + \frac {c_1\tau + c_2}{(1+\tau)^m}
\end{equation}
where $c_1$ and $c_2$ are constants.

Now we take up the problem of completeness of the metrics obtained by Calabi ansatz
starting from a compact $\eta$-Einstein metric. 
First define the function $s(t)$ by
\begin{equation}\label{C27}
s(t) = \int_{\tau_0}^{\tau(t)} \frac{dx}{\sqrt{\varphi(x)}}.
\end{equation}
Then
\begin{equation}\label{C28}
\frac{ds}{dt} = \frac{1}{\sqrt{\varphi(\tau)}}\frac{d\tau}{dt} = \sqrt{\varphi(\tau)}.
\end{equation}
Thus $s(x)$ gives the geodesic length along the $t$-direction with respect to the K\"ahler form $\omega_{\varphi}$ of 
(\ref{C13}); recall $t = \log r$.
\begin{proposition}\label{C29}Let $\omega_{\varphi}$ be the K\"ahler form obtained by
Calabi ansatz starting from a compact Sasaki manifold with 
an $\eta$-Einstein metric $g$. Then $\omega_{\varphi}$ defines a 
complete metric and have noncompact ends towards the end points of $I = (0,b)$
if and only if the following conditions are satisfied at the end points:
\begin{itemize}
\item
At $\tau = 0$, $\varphi$ vanishes at least to the second order.
\item
If $b$ is finite then as $\varphi$ vanishes at $\tau = b$ at least to the second order.
\item
If $b = \infty$ then $\varphi$ grows at most quadratically as $\tau \to \infty$. 
\end{itemize}
\end{proposition}
\begin{proof} First consider at $\tau = 0$. By elementary calculus $s(t) \to \infty$ as $\tau \to 0$
if and only if $\varphi$ vanishes at $0$ at least to the second order. By the same reason, if
$b$ is finite then $\varphi$ must vanish at $\tau = b$ at least to the second order.
Similarly if $b = \infty$,  $s(t) \to \infty$ as $\tau \to \infty$ if and only if $\varphi$ grows at most quadratically.
\end{proof}

\section{Proofs of the theorems}

\begin{proof}[Proof of Theorem \ref{main3}] Going back to (\ref{C26}), assume that $\varphi(0) = 0$
and $\varphi'(0) = 0$. Then $c_1$ and $c_2$ are determined as
\begin{equation}\label{P10}
c_2 = - \frac \kappa{m+1} + \frac c{(m+1)(m+2)}, \quad c_1 = - \kappa + \frac c{m+1}.
\end{equation}
Thus we have
\begin{eqnarray}\label{P11}
\varphi(\tau) &=& \frac{\kappa(1+\tau)}{m+1} + \frac 1{(1+\tau)^m}\left( -\kappa \tau - 
\frac{\kappa}{m+1}\right)
- \frac{c(1+\tau)^2}{(m+1)(m+2)}\\
&& \qquad + \frac c{(1+\tau)^m}\left(\frac\tau{m+1} + \frac1{(m+1)(m+2)}\right).\nonumber
\end{eqnarray}
In the case of this theorem we have $\kappa = 2m + 2$, and we may also assume that
$\kappa$ is any positive number by $D$-homothetic transformation. 
One can check that $\varphi(\tau) > 0$ for all $\tau > 0$ if $c \le 0$, i.e. $b = \infty$, that
$t_1 = -\infty$ for $c \le 0$, $t_2 = \infty$ if $c = 0$ and that $t_2  < \infty$ if $c < 0$.
Using Proposition \ref{C29} one sees that $\omega_{\varphi}$ gives a complete metric.
Let us see when the metric is K\"ahler-Einstein. Recall from (\ref{C18'}) and (\ref{C20}) that
\begin{equation}\label{P2}
\rho_{\varphi} = (\kappa - \frac{m\varphi + (1+\tau)\varphi'}{1+\tau}) \omega^T - 
((\frac{m\varphi}{1+\tau})' + \varphi'')\varphi\,dt\wedge d^ct.
\end{equation}
From this and (\ref{C13}) we see that  $\rho_{\varphi} = \alpha \omega_{\varphi}$
if and only if
\begin{equation}\label{P3}
\kappa - \frac{m\varphi + (1+\tau)\varphi'}{1+\tau} = \alpha(1+\tau),
\end{equation}
\begin{equation}\label{P4}
- (\frac{m\varphi}{1+\tau} + \varphi')' = \alpha.
\end{equation}
Using $\varphi(0) = 0$ and $\varphi'(0) = 0$ one obtains from (\ref{P4})
\begin{equation}\label{P12}
- (\frac{m\varphi}{1+\tau} + \varphi') = \alpha\tau.
\end{equation}
Inserting (\ref{P12}) into (\ref{P3}) we have
\begin{equation}\label{P13}
\kappa = \alpha = \frac c{m+1}.
\end{equation}
But this can not happen because $c \le 0$ and $\kappa > 0$.
The last statement of Theorem \ref{main3} follows from Theorem \ref{s7}. 
This completes the proof of Theorem \ref{main3}.
\end{proof}

\begin{proof}[Proof of Theorem \ref{main5}] Consider the case (a), i.e. $\kappa = 0$. Then
\begin{equation}\label{P14}
\varphi(\tau) =
- \frac{c(1+\tau)^2}{(m+1)(m+2)} + \frac c{(1+\tau)^m}\left(\frac\tau{m+1} + \frac1{(m+1)(m+2)}\right).
\end{equation}
If $c < 0$ then $\varphi(\tau) > 0$ for all $\tau > 0$, i.e. $b = \infty$. Moreover $t_1 = -\infty$
and $t_2 < \infty$. Proposition \ref{C29} shows that $\omega_{\varphi}$ is complete.
One sees from (\ref{P13}) that the metric is not Einstein because $\kappa = 0$ and $c < 0$.

Next consider the case (b), i.e. $\kappa < 0$. Since $\varphi(0) = \varphi'(0) = 0$ and $c < 0$
then (\ref{P11}) shows that there exists a sufficiently large $-c'$ such that for all $c < c'$,
$\varphi_c(\tau) > 0$ for all $\tau > 0$.  Let $c_0$ be the supremum of $c'$ with such a 
property. Let $b$ be the smallest $\tau > 0$ such that $\varphi_{c_0}(\tau) = 0$.
Obviously $\varphi_{c_0}(b) = \varphi'_{c_0}(b) = 0$. Then $t_1 = -\infty$, $t_2 = \infty$
for this $\varphi_{c_0}$. Proposition \ref{C29} shows that this metric is complete.
Now consider $c < c_0$. Then $\varphi_c(\tau) > 0$ for all $\tau > 0$, i.e. $b = \infty$.
One sees that $t_1 = -\infty$, $t_2 < \infty$ and that $\omega_{\varphi}$ is complete
by Proposition \ref{C29}. This metric is Einstein if $c = (m+1)\kappa$ by (\ref{P13}).
This completes the proof of Theorem \ref{main5}.
\end{proof}

\begin{remark}\label{main7} A Sasakian structure is said to be positive (resp. negative)
if $[\rho^T] = \kappa [\omega^T]$ for some positive (resp. negative) real number $\kappa$
 as basic cohomology classes.
A Sasakian structure is said to be null if $[\rho^T] = 0$ as basic cohomology classes.
This case is considered as the case when $\kappa = 0$.
As in the K\"ahler case one can consider in each case the problem of finding a Sasaki metric
with $\rho^T = \kappa \omega^T$, i.e. an $\eta$-Einstein metric, by a transverse K\"ahler deformation. When $\kappa$ is negative
or zero the proofs of the existence results of K\"ahler-Einstein metrics by Aubin and Yau 
apply and one can prove that if the Sasakian structure is negative or null there
exists an $\eta$-Einstein metric with $\rho^T = \kappa \omega^T$ where $ \kappa $ is 
negative or zero. But there are many examples of negative or null Sasaki structures.
For example  the link of an isolated hypersurface singularity defined by a weighted homogeneous polynomials $f$ has a negative Sasakian structure if
$|{\bf w}| - d < 0$ and null Sasakian structure if $|{\bf w}| - d = 0$ where $d$ is the degree of 
the polynomial and ${\bf w} = (w_0, \cdots, w_n)$ is the weights, i.e. 
$$ f(\lambda\cdot\bfz) = f(\lambda^{w_0}z_0, \cdots, \lambda^{w_n}z_n) = \lambda^df(z_0, \cdots,
z_n) = \lambda^d f(\bfz).$$
See \cite{BGM} for more details.
\end{remark}

\section{Construction of Ricci-flat metrics}

In this section we will try to construct a Ricci-flat K\"ahler metric on the total space of the canonical line bundle $K_M$ of
a toric Fano manifold $M$ using by Calabi ansatz. 
But when the Sasaki-Einstein structure is irregular we find it difficult to see if the metric extends smoothly to the zero section. The results in this section are therefore
constructions outside the zero section, just as Proposition \ref{main2}, but when the Sasaki-Einstein structure is regular the metrics obtained extend smoothly to the zero section.

Besides Proposition \ref{main2} we also prove the following results applying the Calabi ansatz to Sasakian $\eta$-Einstein manifolds.
\begin{proposition}\label{main2.5} Let $L$ be a holomorphic line bundle over a compact K\"ahler
manifold $M$ such that $K_M = L^{\otimes p}$ for some positive integer $p$. Let $k$ be a positive
integer and suppose that the total space $S$ of the $U(1)$-bundle associated with $L^{\otimes k}$ satisfies
\begin{equation}\label{main2.6}
\rho^T = \frac{2p}k \omega^T. 
\end{equation}
Then we have the following.
\begin{enumerate}
\item[(a)]
There exists a scalar-flat K\"ahler metric on the total space of $L^{\otimes k}$ minus the zero section, which need not be complete near the zero section 
but is complete near $r = \infty$. This metric is Ricci-flat when $k = p$, that is $L^{\otimes k} = K_M$ and is asymptotic to the
cone metric near infinity. 
\item[(b)]
For any constant $c < 0$ there exists a constant $\gamma > 0$ such that the disk bundle $\{0 < r < \gamma \} \subset L^{\otimes k}$ admits 
a K\"ahler metric of constant scalar curvature $c$. This metric need not be complete near $\{r = 0\}$ but is complete near $\{r = \gamma \}$. 
This metric is K\"ahler-Einstein if $k > p$ and $c = (m+1)(\frac{2p}{k} - 2)$. 
\end{enumerate}
\end{proposition} 

\begin{proof}[Proof of Proposition \ref{main2.5}]\ \ We are in the position of (\ref{C26}) with $\kappa = \frac{2p}k$.
In the case of regular Sasaki-Einstein manifolds, i.e. $U(1)$-bundles associated with the canonical bundle of K\"ahler-Einstein 
manifolds, for $\omega_{\varphi}$ to be complete we need to have $\varphi(0) = 0$ and $\varphi'(0) = 2$. For this fact
refer to for example \cite{Hwang-Singer02}.
We therefore assume the same conditions $\varphi(0) = 0$ and $\varphi'(0) = 2$.  
From these conditions the constants $c_1$ and $c_2$ in (\ref{C26}) are determined
as
\begin{equation}\label{P0}
c_2 = \frac{c}{(m+1)(m+2)} - \frac{\kappa}{m+1},\quad c_1 = \frac c{m+1} + 2 - \kappa.
\end{equation}
Hence $\varphi(\tau)$ is given by
\begin{eqnarray}\label{P1}
\varphi(\tau) &=& \frac{\kappa}{m+1} (1+\tau) - \frac{(\kappa - 2)(m+1)\tau + \kappa}{(1+\tau)^m(m+1)}\\
&& \quad + \frac{c}{(m+1)(m+2)}\left(\frac{(m+2)\tau +1}{(1+\tau)^m} - (1+\tau)^2\right).\nonumber
\end{eqnarray}
One easily checks that if $c \le 0$ then $\varphi(\tau) > 0$ for all $\tau > 0$. 
Thus $b = \infty$, i.e. $I = (0,\infty)$. Since $\varphi(0) = 0$ and
$\varphi'(0) = 2$ we see from (\ref{C8}) that $t_1 = - \infty$. On the other hand if $c = 0$ then
$\varphi(\tau)$ has linear growth as $\tau \to \infty$. Thus $t_2 = \infty$. If $c < 0$ then 
$\varphi(\tau)$ has quadratic growth as $\tau \to \infty$. Thus $t_2 < \infty $. In either case
$\omega_{\varphi}$ defnes a complete metric near the infinity since the growth is at most
quadratic.

In the case of $c = 0$, $\omega_{\varphi}$ defines a scalar-flat K\"ahler metric on
the $L$ minus the zero section since $t_2 = \infty$, and the metric is complete near the infinity. 
In the case of $c < 0$, $\omega_{\varphi}$ defines a 
K\"ahler metric of constant negative scalar curvature $c$ on a disk bundle of $L$ minus the zero section since $t_2 < \infty$.

This metric satisfies $\rho_{\varphi} = \alpha \omega_{\varphi}$ if and only if (\ref{P3}) and (\ref{P4}) hold.
From (\ref{P4}), $\varphi(0) = 0$ and $\varphi'(0) = 2$ we see
\begin{equation}\label{P5}
- (\frac{m\varphi}{1+\tau} + \varphi') = \alpha\tau - 2.
\end{equation}
Inserting (\ref{P5}) into (\ref{P3}) and using $\kappa = \frac{2p}k$ we have
\begin{equation}\label{P6}
\alpha = \frac{2p}k -2.
\end{equation}
But since the scalar curvature is computed as
\begin{equation}\label{P7}
c = (m+1)\alpha
\end{equation}
and since $c \le 0$ then
the metric is Einstein if it happens that
$k \ge p$ and that
\begin{equation}\label{P8}
c = (m+1)(\frac{2p}{k} - 2).
\end{equation}
If $c = 0$ this certainly happens when $k = p$.
This completes the proof of Theorem \ref{main2.5}.
\end{proof}

\begin{proof}[Proof of Proposition \ref{main2}]\ \ By Theorem \ref{s7}, $S$ admits an
$\eta$-Einstein Sasaki metric with $\rho^T = (2p/k) \omega^T$. Hence we can apply 
Theorem \ref{main2.5}. This completes the proof.
\end{proof}

\begin{proposition}\label{main6}Let $L$ be a holomorphic line bundle over a compact K\"ahler
manifold $M$ such that $K_M = L^{\otimes p}$ for some positive integer $p$. Let $k$ be a positive
integer and suppose that
the total space $S$ of the $U(1)$-bundle associated with $L^{\otimes k}$ has a
Sasakian $\eta$-Einstein metric with
\begin{equation}\label{main2.6'}
\rho^T = \kappa \omega^T
\end{equation}
for some non-positive constant $\kappa$. Then we have the following.\\
\noindent
Case $\kappa = 0$\ :\ \ 
There exists a scalar-flat
K\"ahler metric on the total space of $L^{\otimes k}$ minus the zero section. 
Further,
for any constant $c < 0$ there exists a constant $\gamma > 0$ such that the disk bundle
$\{0< r < \gamma \} \subset L^{\otimes k}$ admits a K\"ahler metric of constant scalar curvature $c$.
This metric need not be complete along the zero section and complete near $r = \gamma$.
This metric is K\"ahler-Einstein if $c = -2(m+1)$. \\
\noindent
Case $\kappa < 0$\ :\ \ 
We have a negative constant $c_0$ such that there exists a 
K\"ahler metric of negative constant scalar curvature on the total space of $L^{\otimes k}$ minus the zero section
and that for any $c < c_0$ there exists a constant $\gamma > 0$ such that the disk bundle
$\{ 0 < r < \gamma \} \subset L^{\otimes k}$ admits a K\"ahler metric of constant scalar 
curvature $c$. This metric need not be complete near the zero section and complete near $r = \gamma$. 
This metric is K\"ahler-Einstein if $k > p$ and $c = (m+1)(\kappa - 2)$. 
\end{proposition} 
\begin{proof}
Consider the case $\kappa = 0$. By (\ref{P1}) we have
\begin{equation}\label{P9}
\varphi(\tau) = \frac{ 2\tau}{(1+\tau)^m}
 + \frac{c}{(m+1)(m+2)}\left(\frac{(m+2)\tau +1}{(1+\tau)^m} - (1+\tau)^2\right).
\end{equation}
As in the proof of Proposition \ref{main2.5} if $c \le 0$ we have $b = \infty$ and $t_1 = - \infty$.
If $c = 0$ we have $t_2 = \infty$, and if $c < 0$ we have $t_2 < \infty$. The rest of the proof
goes as in the proof of Proposition \ref{main2.5}. The metric
is Einstein if and only if $c = -2(m+1)$. \\
Consider the case $\kappa < 0$. 
Let $\varphi_c$ denote $\varphi$ for a given $c $. (\ref{P1}) shows that,
since $\varphi(0) = 0$ and $\varphi'(0) = 2$ then there exists a sufficiently
large $-c'$ such that for any $c < c'$ we have $\varphi_c(\tau) > 0$ for all $\tau > 0$. 
Let $c_0$ be the supremum
of $c'$'s with such a property. 
Let $b$ be the smallest $\tau>0$ such that $\varphi_{c_0}(\tau) = 0$. Obviously $\varphi_{c_0}'(b) = 0$.
Thus, as in Proposition \ref{C29}, $\omega_{\varphi_0}$ defines a complete metric. In this case $t_1 = -\infty$
and $t_2 = \infty$, so the metric is defined on the whole total space of $L^{\otimes k}$ minus the zero section. 
Take $c$ such that  $c < c_0$, then $\varphi_c > 0$ for all $\tau > 0$, so $b = \infty$.
Since
$\varphi_c$ grows quadratically as $\tau \to \infty$ then $\omega_{\varphi_c}$ is complete by 
Proposition \ref{C29} except near the zero section. In this case $t_2 < \infty$. The metric
is Einstein if and only if $c = (m+1)(\kappa -2)$. That this metric is asymptotic to the cone metric will be proved in 
Lemma \ref{asymp} below. This completes the proof
\end{proof}

In the first version of arXiv:math/0703138 posted under a different title, the author claimed that the metric obtained in 
Proposition \ref{main2} extends smoothly to the zero section. But its proof is not correct because of the following reasons.
First of all it is proved in \cite{FOW} that if the Reeb vector field corresponds to $\xi \in \mathfrak g$ then the K\"ahler potential
$F_{\xi}^{can}$ is given by the equation (61) in  \cite{FOW}:
$$ F_{\xi}^{can} = \frac{r^2}2 = \frac12 l_{\xi}(y) $$
where $y$ is the position vector on the moment cone. Here $y$ is considered as a part of the action-angle coordinates
with fixed symplectic form and varying complex structure. If we choose another Reeb vector field corresponding to $\xi^{\prime} \in \mathfrak g$
then the new K\"ahler potential is given by
\begin{equation}\label{s11'}
\frac{r^{\prime2}}2 = \frac12 r^2 \frac {l_{\xi'}(y)}{l_{\xi}(y)} = \frac12 r^2 \frac {\la \xi',y \ra}{\la\xi,y \ra}.
\end{equation}
But in the Abreu-Guillemin arguments when the K\"ahler potential varies the complex structure varies, which means that
the holomorphic coordinates change and also $i\p \bp$-operator changes. As a result, to recover the new K\"ahler form,
we have to apply new $i\p^{\prime} \bp^{\prime}$ to the new K\"ahler potential 
$\frac{r^{\prime2}}2$. This implies that in the logarythmic affine coordinates with fixed complex structure and varying symplectic structure
the two K\"ahler forms with varying Reeb vector fields are not related by 
$$ r'^2 = r^2 \exp \psi$$
with $\psi \in C^{\infty}(S)$, although $\frac{r^{\prime2}}2$ and $\frac12 r^2$ are related in the action-angle coordinates by $\frac {\la \xi',y \ra}{\la\xi,y \ra} \in C^{\infty}(S)$
as in (\ref{s11'}).
Second of all if two Sasaki structures are related in the standard holomorphic coordinates (with fixed complex structure and varying symplectic structure) by
$ r^{\prime} = r\exp\vph$
with a basic function $\vph \in C^{\infty}(S)$ then the two Sasaki structures share a common Reeb vector field. When we consider toric Sasaki manifolds all metrics are
$T^{m+1}$-invariant and if two Sasaki structures are related by $ r'^2 = r^2 \exp \psi$
with $\psi \in C^{\infty}(S)$ then $\psi$ is necessarily basic. This implies that if two Sasaki structures have different Reeb vector fields then the two Sasaki structures are not
related by $ r'^2 = r^2 \exp \psi$
with $\psi \in C^{\infty}(S)$. Therefore the when the Futaki invariant of $M$ is not zero the Sasaki-Einstein structure on $S$ is never bundle-adapted in the sense of the version 1
of arXiv:math/0703138.

\begin{lemma}\label{asymp} The Ricci-flat metric in (a) of Proposition \ref{main2.5} is asymptotic to the cone metric near $r = \infty$.
\end{lemma}
\begin{proof}

Putting in (\ref{P1}) $\kappa = 2$, $c = 0$ we have 
$$
\varphi(\tau) = \frac 2{m+1} ((1+ \tau) - \frac 1{(1+\tau)^m}) = \frac 2{m+1} \frac{(1+ \tau)^{m+1} -1}{ (1+\tau)^m}.
$$
Take $\tau_0 = 2^{\frac 1{m+1}} -1$. Then $\tau(t)$ is obtained from (\ref{C8}) as
\begin{eqnarray}
t &=& \int_{\tau_0}^{\tau(t)} \frac{dt}{\varphi(\tau)} = \frac{m+1}2 \int_{\tau_0}^{\tau(t)} \frac{(\tau + 1)^m}{(\tau+1)^{m+1} -1} dt \\
&=& \frac 12 \log ((\mu(t) + 1)^{m+1} -1.
\end{eqnarray}
Since $t = \log r$ then we have 
$$ e^{2t} = r^2 = (\mu(t) + 1)^{m+1} - 1,$$
$$\tau(t) = (e^{2t} + 1)^{\frac 1{m+1}} - 1.$$
$F(t)$ is obtained from (\ref{C9}) as 
$$
F(t) = \int_{\tau_0}^{\mu(t)} \frac {m+1}2 \frac {(\tau +1)^m\tau}{(\tau+1)^{m+1} - 1} d\tau.
$$
Hence we have 
\begin{eqnarray}
\tau &=& F'(t) = \frac {m+1}2 \frac {(\tau(t) +1)^m\tau(t)}{(\tau(t)+1)^{m+1} - 1} \tau'(t)\\
&=& (e^{2t} + 1)^{\frac 1{m+1}} - 1.
\end{eqnarray}

Let $\zeta$ be the $(m+1)$-root of unity. Then 
\begin{eqnarray}
 \sum_{j=0}^m \zeta^j \Pi_{i\ne j} (x - \zeta^i) &=& - \sum_{j=0}^m (x- \zeta^j) \Pi_{i\ne j} (x - \zeta^i) + \sum_{j=0}^m x\Pi_{i\ne j}(x - \zeta^i) \nonumber\\
 &=& - (m+1)(x^{m+1} -1) + x (x^{m+1} - 1)' \nonumber \\
 &=& -(m+1))(x^{m+1} -1) + (m+1)x^{m+1} \nonumber\\
 &=& m+1.\nonumber
\end{eqnarray}
Using this we get 
\begin{eqnarray}
&&\frac d{dt} ((e^{2t} + 1)^{\frac 1{m+1}} + \frac 1{m+1} \sum_{j=0}^m \zeta^j \log((e^{2t} + 1)^{\frac 1{m+1}} - \zeta^j))\nonumber\\
&=& \frac 2{m+1} e^{2t}(e^{2t }+1)^{- \frac m{m+1}}(1 + \frac 1{m+1} \sum_{j=0}^m \frac{\zeta^j}{(e^{2t} + 1)^{\frac 1{m+1}} - \zeta^j})\nonumber\\
&=& \frac 2{m+1}  e^{2t}(e^{2t} +1)^{- \frac m{m+1}}(1 + \frac 1{m+1} \frac 1{e^{2t}} (m+1))\nonumber\\
&=& \frac 2{m+1}(1 + e^{2t})^{\frac 1{m+1}} = \frac 2{m+1} (F'(t) + 1).\nonumber
\end{eqnarray}
Thus K\"ahler potential is expressed up to the mutiple of $\frac {m+1}2$ as
$$ (e^{2t} + 1)^{\frac 1{m+1}} + \frac 1{m+1} \sum_{j=0}^m \zeta^j \log((e^{2t} + 1)^{\frac 1{m+1}} - \zeta^j). $$
Since we are using the $\eta$-Einstein metric such that $\mathrm{Ric}^T = 2 \omega^T$, to deform it into the Sasaki-Einstein metric
such that $\mathrm{Ric}^{\prime T} = 2(m+1) \omega^{\prime T}$ we need to perform $D$-homothetic transformation with 
$$ r = \tilde r^{m+1}.$$
Thus the K\"ahler potential of our Ricci-flat K\"ahler metric is given by
$$ f = ((\tilde r^{2m+2} + 1)^{\frac 1{m+1}} + \frac 1{m+1} \sum_{j=0}^m \zeta^j \log((\tilde r^{2m+2} + 1)^{\frac 1{m+1}} - \zeta^j)). $$
Putting $s = 1/\tilde r$ and using l'H\^optal's theorem we have
\begin{eqnarray}
\lim_{\tilde r \to \infty} \frac {f - \tilde r^2}{\tilde r^{\alpha}} &=& \lim_{\tilde r \to \infty} \frac{2(\tilde r^{2m+2} +1)^{\frac 1{m+1}} - 2\tilde r^2}{\alpha \tilde r^{\alpha}}\nonumber\\
&=&\lim_{s \to 0}  \frac{2(1 +s^{2m+2} )^{\frac 1{m+1}} - 2}{\alpha  s^{2-\alpha}}.\nonumber
\end{eqnarray}
Using
$$
(1 + T)^{\frac 1{m+1}} = 1 + \frac 1{m+1} T - \frac m{2(m+1)^2} T^2 + ...
$$
further the limit above converges when $\alpha = -2m$ as
$$ \lim_{s \to 0} \frac{\frac 2{m+1} s^{2m+2} - ...}{\alpha s^{2-\alpha}} = - \frac 1{m(m+1)}.$$
We finally get
$$ f = \tilde r^2 - \frac 1{m(m+1)}\tilde r^{-2m} + O(\tilde r^{-4m-2}).$$
This completes the proof of Lemma.
\end{proof}

\end{document}